\documentclass{amsart}
\usepackage{amscd,amssymb}
\usepackage{latexsym}

\newtheorem{theorem}{Theorem}[section]
\newtheorem{proposition}[theorem]{Proposition}
\newtheorem{lemma}[theorem]{Lemma}
\newtheorem{corollary}[theorem]{Corollary}

\theoremstyle{definition}
\newtheorem{definition}[theorem]{Definition}

\theoremstyle{remark}

\numberwithin{equation}{section}

\DeclareMathOperator{\colim}{colim}

\DeclareMathOperator{\Hom}{Hom}
\DeclareMathOperator{\Ext}{Ext}
\DeclareMathOperator{\ho}{Ho}
\DeclareMathOperator{\Ch}{Ch}
\DeclareMathOperator{\Spec}{Spec}

\DeclareMathOperator{\QCo}{QCo}
\DeclareMathOperator{\ann}{ann}
\DeclareMathOperator{\cofinal}{cofin}

\newcommand{\qch}{\Ch _{\QCo }}

\newcommand{\cat}[1]{\mathcal{#1}}
\newcommand{\Mod}{\text{-Mod}}
\newcommand{\Sh}{\cat{S}h}
\newcommand{\Z}{\mathbb{Z}}

\newcommand{\ideal}[1]{\mathfrak{#1}}

\newcommand{\mathcolon}{\colon\,}
\newcommand{\ie}{\textit{i.e.}}

\newcommand{\proj}{\text{-proj}}
\newcommand{\cof}{\text{-cof}}
 
\newcommand{\inj}{\text{-inj}}

\newcommand{\boxprod}{\Box }
\newcommand{\rlp}{right lifting property with respect to }
\newcommand{\llp}{left lifting property with respect to }
\newcommand{\fgd}{finite global dimension}
\newcommand{\fhgd}{finite hereditary global dimension}
\begin{document}

\title[Chain complexes of sheaves]{Model category structures on chain
complexes of sheaves}

\date{\today}

\author{Mark Hovey}
\address{Department of Mathematics \\ Wesleyan University
\\ Middletown, CT 06459}
\email{hovey@member.ams.org}


\begin{abstract}
In this paper, we try to determine when the derived category of an
abelian category is the homotopy category of a model structure on the
category of chain complexes.  We prove that this is always the case when
the abelian category is a Grothendieck category, as has also been done
by Morel.  But this model structure is not very useful for defining
derived tensor products.  We therefore consider another method for
constructing a model structure, and apply it to the category of sheaves
on a well-behaved ringed space.  The resulting flat model structure is
compatible with the tensor product and all homomorphisms of ringed
spaces.  
\end{abstract}

\maketitle

\section*{Introduction}\label{sec-intro}

It very often happens in mathematics that one has a category $\cat{C}$
and a collection of maps $\cat{W}$ in $\cat{C}$ that one would like to
consider as isomorphisms.  In this situation, one can formally invert
the maps in $\cat{W}$, but the resulting localization $\ho \cat{C}$ of
$\cat{C}$ may not be a category in general, because $\ho \cat{C}(X,Y)$
may not be a set.  Furthermore, it is hard to get a handle on maps in
$\ho \cat{C}$ from $X$ to $Y$.  Model categories were invented by
Quillen~\cite{quillen-htpy} to get around these problems.  In general,
it is hard to prove a given category $\cat{C}$ is a model category, but,
having done so, many structural results about $\ho \cat{C}$ follow
easily; for example, $\ho \cat{C}$ is then canonically enriched over the
homotopy category of simplicial sets.  And of course one can then use
the considerable body of results about model categories to investigate
$\cat{C}$.

An obvious example of a situation where one wants to invert some maps is
the construction of the derived category of an abelian category
$\cat{A}$.  Recall that this is the localization of the category $\Ch
(\cat{A})$ of (unbounded) chain complexes by maps which induce homology
isomorphisms.  The category of nonnegatively graded chain complexes of
$R$-modules was one of Quillen's first examples of a model category.
Nevertheless, the first published proof that the category of unbounded
chain complexes of $R$-modules is a model category appears to be
in~\cite{hovey-model}.  

We begin the paper by establishing a model structure on $\Ch (\cat{A})$
whose homotopy category is the derived category, when $\cat{A}$ is a a
Grothendieck category.  In particular, $\cat{A}$ could be the category
of sheaves on a ringed space, or the category of quasi-coherent sheaves
on a quasi-compact and quasi-separared scheme.  The author has heard
that this has been done by F. Morel as well, but does not know any
details.  The injective model structure is natural for exact functors of
abelian categories, but not for left adjoints which are only right
exact.  Also, if $\cat{A}$ is closed symmetric monoidal, the injective
model structure will not be compatible with the tensor product, making
it of no use for defining the derived tensor product.

We therefore discuss a different method for constructing a model
structure on $\Ch (\cat{A})$.  This method enables us to define a
different model structure on $\Ch (\cat{A})$ in case $\cat{A}$ has a set
of generators of finite projective dimension.  In particular, we apply
it when $\cat{A}$ is the category of quasi-coherent sheaves on a nice
enough scheme, using the locally free sheaves as the generators.  Though
the resulting locally free model structure is still not compatible with
the tensor product, it does give us some information about the resulting
derived category that does not seem accessible from the injective model
structure.  

But this method works better when $\cat{A}$ is the category of sheaves
on a ringed space $(S,\cat{O})$ satisfying a hypothesis related to
finite cohomological dimension.  In this case, we construct a flat model
structure on $\Ch (\cat{A})$ that is compatible with the tensor product.
We then get model categories of differential graded $\cat{O}$-algebras
and of differential graded modules over a given differential graded
$\cat{O}$-algebra.  The flat model structure is also natural for
arbitrary maps of ringed spaces.  

To understand this paper, the reader needs to know some basic facts
about model categories, Grothendieck categories, and sheaves.  A good
introduction to model categories is~\cite{dwyer-spalinski}.  The
book~\cite{hovey-model} is a more in-depth study, but still starting
from scratch.  All the terms we need are defined in~\cite{hovey-model};
we will give specific references as needed.  For Grothendieck
categories, \cite{stenstrom} is sufficient.  For sheaves and schemes, we
try to refer mostly to~\cite{hartshorne}, but we also need more advanced
results occasionally.  

The author would like to thank Matthew Ando and Amnon Neeman for helpful
discussions about sheaves, and Dan Christensen for his many useful
suggestions.

\section{Grothendieck categories}\label{sec-grot} 

In this section, we develop the basic structural properties of
Grothendieck categories that we need.  In particular, we show that every
object in a Grothendieck category is small and that the fundamental
lemma of homological algebra holds in a Grothendieck category.  

Recall that a cocomplete abelian category is called an AB5 category if
directed colimits, or, equivalently, filtered colimits, are exact.  A
Grothendieck category is an AB5 category that has a generator.  Recall
that $U$ is a generator of $\cat{A}$ if the functor $\cat{A}(U,-)$ is
faithful.  For example, the category of sheaves on a ringed space, or a
ringed topos, is a Grothendieck category.  The Popescu-Gabriel theorem,
a proof of which can be found in~\cite[Section~X.4]{stenstrom}, asserts
that every Grothendieck category $\cat{A}$ is equivalent to the
full subcategory of $\cat{T}$-local objects in a module category $R\Mod $,
for some hereditary torsion theory $\cat{T}$ and some ring $R$.  The
ring $R$ can be taken to be $\Hom (U,U)$ for a generator $U$, but since
the generator $U$ is not canonically attached to $\cat{A}$, neither is
the ring $R$.  

One of the most basic tools used to establish the existence of a model
category structure is the small object
argument~\cite[Section~2.1]{hovey-model}.  In order to use the small
object argument, we need to know that the objects in a Grothendieck
category are small; \ie that maps out of an object commute with long
enough colimits.  In order to make this precise, we need some
definitions.

\begin{definition}\label{defn-small}
\begin{enumerate}
\item Given a limit ordinal $\lambda $, the \emph{cofinality} of
$\lambda $, $\cofinal \lambda $, is the smallest cardinal $\kappa $ such
that there exists a subset $T$ of $\lambda $ with $|T|=\kappa $ and
$\sup T=\lambda $.
\item Given an object $A$ in a cocomplete category $\cat{C}$ and a cardinal
$\kappa $, we say that $A$ is \emph{$\kappa $-small} if, for every
ordinal $\lambda $ with $\cofinal \lambda > \kappa $ and every
colimit-preserving functor $X\mathcolon \lambda \xrightarrow{}\cat{C}$,
the natural map $\colim_{i <\lambda }
\cat{C}(A,X_{i})\xrightarrow{}\cat{C}(A,\colim _{i<\lambda }X_{i})$ is
an isomorphism.  
\item An object $A$ in a cocomplete category $\cat{C}$ is called
\emph{small} if it is $\kappa $-small for some cardinal $\kappa $.  
\end{enumerate}
\end{definition}

\begin{proposition}\label{prop-grot-small}
Every object in a Grothendieck category is small.  
\end{proposition}

We do not know if this proposition holds more generally,
\texttt{e.g.}, in any AB5 category, but it seems unlikely.  

\begin{proof}
We may as well assume that our Grothendieck category $\cat{A}$ is the
localization of $R\Mod $ with respect to a hereditary torsion theory
$\cat{T}$, for some ring $R$.  Let $\kappa $ be the larger of $\infty $
and the cardinality of $R$, let $\lambda $ be an
ordinal with $\cofinal \lambda > \kappa $, and let $X\mathcolon
\lambda \xrightarrow{}\cat{A}$ be a colimit-preserving functor.  We will
first show that $\colim X_{i}$, calculated in $R\Mod $, is still
$\cat{T}$-local, so is also the colimit in $\cat{A}$.  This proof will
depend on the fact that both $R/\ideal{a}$ and $\ideal{a}$ are $\kappa
$-small in $R\Mod $~\cite[Example~2.1.6]{hovey-model}, for all (left)
ideals $\ideal{a}$ of $R$.

To see this, first note that $\colim X_{i}$ is torsion-free. Indeed,
$\cat{T}$ is generated by cyclic modules $R/\ideal{a}$, so it suffices
to show that $R\Mod (R/\ideal{a},\colim X_{i})=0$.  But we have chosen
$\kappa $ so that 
\[
R\Mod (R/\ideal{a},\colim X_{i}) \cong \colim R\Mod (R/\ideal{a},X_{i})
=0,
\]
since each $X_{i}$ is torsion-free.  Hence $\colim X_{i}$ is
torsion-free. 

It follows that the localization of $\colim X_{i}$ is 
\[
\colim _{\ideal{a}} \Hom (\ideal{a},\colim_{i} X_{i}),
\]
where the colimit is taken over ideals $\ideal{a}$ such that
$R/\ideal{a}$ is torsion, as in~\cite[Section~IX.1]{stenstrom}.  But
then we have
\begin{align*}
\colim _{\ideal{a}} \Hom (\ideal{a},\colim _{i}X_{i}) & \cong  
\colim _{\ideal{a}} \colim _{i} \Hom (\ideal{a},X_{i})  \\
& \cong \colim _{i}\colim _{\ideal{a}} \Hom (\ideal{a},X_{i}) = 
\colim _{i}X_{i}.  
\end{align*}
Thus $\colim X_{i}$ is already local.  

Now suppose $M$ is an arbitrary local module.  There is a cardinal
$\kappa '$ such that $M$ is $\kappa '$-small as an $R$-module, and we
can choose $\kappa '\geq \kappa $.  It is then immediate from the
argument above that $M$ is $\kappa '$-small in $\cat{A}$.  
\end{proof}

We will be most interested in the category of unbounded chain complexes
$\Ch (\cat{A})$ on a Grothendieck category $\cat{A}$.  This is again a
Grothendieck category; since colimits are taken dimensionwise, filtered
colimits are obviously exact, and if $U$ is a generator of $\cat{A}$,
then the disks $D^{n}U$ are generators of $\Ch (\cat{A})$.  Recall that
$D^{n}U$ is the complex which is $U$ in degrees $n$ and $n-1$ and $0$
elsewhere, with the interesting differential being the identity map.
To see that the disks generate $\Ch (\cat{A})$, use the adjunction
relation $\Ch (\cat{A})(D^{n}U,X)\cong \cat{A}(U,X_{n})$.  In
particular, every object of $\Ch (\cat{A})$ is small.

We recall the fundamental lemma of homological algebra, which holds in
any abelian category~\cite[Section~VI.8]{mitchell-cat}.  

\begin{lemma}\label{lem-grot-fund}
Suppose $\cat{A}$ is an abelian category, and
\[
0\xrightarrow{} A \xrightarrow{f}B \xrightarrow{g} C \xrightarrow{} 0
\]
is a short exact sequence in $\Ch (\cat{A})$.  Then there is a natural
long exact sequence 
\[
\dots \xrightarrow{} H_{n+1}C \xrightarrow{\partial _{*}} H_{n}A
\xrightarrow{f_{*}} H_{n}B \xrightarrow{g_{*}} H_{n}C
\xrightarrow{\partial _{*}} H_{n-1}A \xrightarrow{} \dots
\]
in homology.  
\end{lemma}

Recall that a map $f$ in $\Ch (\cat{A})$ is called a
\emph{quasi-isomorphism} if $H_{*}f$ is an isomorphism.  The following
corollary is immediate.  

\begin{corollary}\label{cor-proper}
Suppose $\cat{A}$ is an abelian category, and
\[
\begin{CD}
A @>f>> B \\
@VrVV @VVsV \\
C @>>g> D
\end{CD}
\]
is a commutative square in $\Ch (\cat{A})$.  
\begin{enumerate}
\item [(a)] If the square above is a pushout square and $f$ is an
injective quasi-isomorphism, so is $g$.  
\item [(b)] If the square above is a pushout square, $r$ is injective
and $f$ is a quasi-isomorphism, then $g$ is a quasi-isomorphism.
\item [(c)] If the square above is a pullback square and $g$ is a
surjective quasi-isomorphism, so is $f$.  
\item [(d)] If the square above is a pullback square, $s$ is surjective,
and $g$ is a quasi-isomorphism, then $f$ is a quasi-isomorphism.  
\end{enumerate}
\end{corollary}

It is also useful to know that homology commutes with colimits. 

\begin{lemma}\label{lem-colimits-homology}
Suppose $\cat{A}$ is an AB5 category, $\cat{I}$ is a small filtered
category, and $F\mathcolon \cat{I}\xrightarrow{}\Ch (\cat{A})$ is a
functor.  Then there is a natural isomorphism $\colim H_{n}F
\xrightarrow{} H_{n}\colim F$.
\end{lemma}

\begin{proof}
We have an exact sequence 
\[
0 \xrightarrow{} Z_{n}F(i) \xrightarrow{}
F(i)_{n}\xrightarrow{d}F(i)_{n-1} \xrightarrow{} B_{n}F(i) \xrightarrow{} 0
\]
for each object $i$ of $\cat{I}$.  Since filtered colimits are exact, we
find that $Z_{n}\colim _{i} F(i)\cong \colim _{i}Z_{n}F(i)$ and
smiilarly for $B_{n}$.  Applying colimits to the short exact sequences 
\[
0\xrightarrow{}B_{n}F(i) \xrightarrow{} Z_{n}F(i) \xrightarrow{}H_{n}F(i)
\]
completes the proof.  
\end{proof}

In particular, using transfinite induction, we obtain the following
proposition.  

\begin{proposition}\label{prop-colimits-weak-equivs}
Suppose $\cat{A}$ is an AB5 category, $\lambda $ is an ordinal,
and $X\mathcolon \lambda \xrightarrow{}\Ch (\cat{A})$ is a
colimit-preserving functor such that, for all $\alpha < \lambda $, the
map $X_{\alpha }\xrightarrow{}X_{\alpha +1}$ is a quasi-isomorphism.
Then the map $X_{0}\xrightarrow{}\colim _{\alpha <\lambda }X_{\alpha }$
is a quasi-isomorphism.
\end{proposition}

In the theory of model categories, one frequently has a set of maps $J$
and wants to consider the class $J\inj $ of maps that look like
fibrations to $J$ and the class $J\cof $ of maps that look like
cofibrations to $J\inj $.  These classes are defined by lifting
properties~\cite[Section~2.1]{hovey-model}.  

\begin{corollary}\label{cor-homology-J}
Let $\cat{A}$ be a Grothendieck category.  Suppose $J$ is a set of
injective quasi-isomorphisms in $\Ch (\cat{A})$.  Then $J\cof $ consists
of injective quasi-isomorphisms.
\end{corollary}

\begin{proof}
By the small object argument~\cite[Theorem~2.1.14 and
Corollary~2.1.15]{hovey-model}, every element of $J\cof $ is a retract
of a transfinite composition of pushouts of elements of $J$.  Injections
are closed under retracts and pushouts in any abelian category; the AB5
condition guarantees that they are also closed under transfinite
compositions.  Part~(a) of Corollary~\ref{cor-proper} then shows that
pushouts of maps of $J$ are quasi-isomorphisms.
Proposition~\ref{prop-colimits-weak-equivs} shows that transfinite
compositions of quasi-isomorphisms are quasi-isomorphisms.  It is clear
that retracts of injective quasi-isomorphisms are quasi-isomorphisms.
\end{proof}

Note that this corollary holds in any AB5 category as long as the
domains and codomains of the maps of $J$ are small, so that the small
object argument applies.  

\section{The injective model structure}\label{sec-inj}

In this section, we construct the injective model structure on $\Ch
(\cat{A})$ when $\cat{A}$ is a Grothendieck category.

\begin{definition}\label{defn-grot-fib}
Define a map $p\mathcolon X\xrightarrow{}Y$ in $\Ch (\cat{A})$ to be an
\emph{injective fibration} if it has the \rlp all injective weak
equivalences in $\Ch (\cat{A})$.  
\end{definition}

Note that, by definition, a complex $X$ is injectively fibrant if and
only if $X$ is $DG$-injective in the sense
of~\cite[Section~7]{avramov-foxby-halperin}.  The arguments in that
paper then show that $X$ is $DG$-injective if and only if each $X_{n}$
is injective in $\cat{A}$ and $X$ is $K$-injective in the sense of
Spaltenstein~\cite{spaltenstein}.  We will see in
Proposition~\ref{prop-inj-fibs} that an injective fibration is just a
dimensionwise split surjection with $DG$-injective kernel.

Then the object of this section is to prove the following theorem.

\begin{theorem}\label{thm-inj}
Suppose $\cat{A}$ is a Grothendieck category.  Then the derived
category of $\cat{A}$ is the homotopy category of a cofibrantly
generated proper model structure on $\Ch (\cat{A})$ where the
cofibrations are the injections, the fibrations are the injective
fibrations, and the weak equivalences are the quasi-isomorphisms.
\end{theorem}

We call this model structure on $\Ch (\cat{A})$ the \emph{injective
structure}, or the \emph{injective model structure}.  


\begin{corollary}\label{cor-grot-sheaves}
\begin{enumerate}
\item [(a)] Suppose $(S,\cat{O})$ is a ringed space \textup{(}or a
ringed topos\textup{)}.  Then the injective structure on the category
$\Ch (\cat{O}\Mod )$ of unbounded complexes of sheaves of
$\cat{O}$-modules is a model structure, whose homotopy category is the
derived category.
\item [(b)] Suppose $S$ is a quasi-compact, quasi-separated scheme.
Then the injective structure on the category $\qch (S)$ of unbounded
complexes of quasi-coherent sheaves of $\cat{O}_{S}$-modules is a model
structure, whose homotopy category is the derived category of
quasi-coherent sheaves.  
\end{enumerate}
\end{corollary}

\begin{proof}
It is well-known that the category of sheaves on a ringed space is a
Groth\-en\-dieck category~\cite[Proposition~3.1.1]{grothendieck-tohoku}.
The category of quasi-coherent sheaves on ringed space is an abelian
subcategory of all sheaves, closed under colimits.  So colimits are
exact.  It remains to show that the category of quasi-coherent sheaves
on a quasi-compact, quasi-separated scheme has a generator.  This is a
corollary of Deligne's result in~\cite[Appendix,
Prop.~2]{hartshorne-residues}, which asserts that every quasi-coherent
sheaf is the colimit of finitely presented sheaves.  Since there is only
a set of finitely presented sheaves, the direct sum of all of them will
serve as a generator.
\end{proof}

Theorem~\ref{thm-inj} is of course a generalization of the corresponding
theorem about complexes of modules~\cite[Theorem~2.3.13]{hovey-model},
and will be proved in a similar way.  We need sets $I$ and $J$ of
generating cofibrations and generating trivial cofibrations.  We will
find these sets by just taking all injections (resp. all injective weak
equivalences) whose cardinality is not too large.  Then we will show
that $I\cof $ is the class of injections and $J\cof $ is the class of
injective weak equivalences.  We must also show that every map in $I\inj
$ is an injective fibration and a weak equivalence.  The recognition
theorem~\cite[Theorem~2.1.19]{hovey-model} will then prove that we do
get a model category.  Properness is automatic from
Corollary~\ref{cor-proper}.

To carry out this plan, we need a notion of cardinality.  This notion
will depend on a way of representing our Grothendieck category $\cat{A}$
as a localization of a category of modules, but the resulting model
structure will be independent of this choice.  So throughout the rest of
this section, we will assume that $\cat{A}$ is the localization of
$R\Mod $, for some ring $R$, with respect to a hereditary torsion
theory.  In particular, we will think of objects of $\cat{A}$ as being
$R$-modules.  

\begin{definition}\label{defn-card}
Suppose $\cat{A}$ is a Grothendieck category.  Define the
\emph{cardinality} of $X$, $|X|$, to be the cardinality of $X$ as an
$R$-module.  Given a chain complex $X\in \Ch (\cat{A})$, we define $|X|$
to be the cardinality of the disjoint union of the $X_{n}$.  Define
$\gamma $ to be the supremum of $\infty $ and $2^{|R|}$.  Then define
$I$ to be a set containing one element of each isomorphism class of
injections $A\xrightarrow{}B$ in $\Ch (\cat{A})$ with $|B|\leq \gamma $.
Define $J$ to be the set of all quasi-isomorphisms in $I$.
\end{definition}

The reason for choosing $\gamma $ as we have done is the following
lemma.  

\begin{lemma}\label{lem-cardinality}
Suppose $R$ is a ring, $\cat{T}$ is a torsion theory on $R\Mod $, and
\[
L\mathcolon R\Mod \xrightarrow{}R\Mod 
\]
is the corresponding localization functor.  If $M$ is an $R$-module with
$|M|\leq \gamma $, then $|LM|\leq \gamma $.
\end{lemma}

\begin{proof}
We may as well assume $M$ is torsion-free, since killing the torsion
only decreases the cardinality of $M$, without changing $LM$.  Then
$LM=\colim \Hom (\ideal{a},M)$, where the colimit is taken over ideals
$\ideal{a}$ such that $R/\ideal{a}$ lies in $\cat{T}$.  Thus 
\[
|LM|\leq 2^{|R|} |M|^{|R|}. 
\]
In case $|R|$ is finite, $\gamma =\infty $, and one can easily see from
this equation that $|LM|$ is countable when $|M|$ is so.  In case $R$ is
infinite, we have 
\[
|LM|\leq 2^{|R|} (2^{|R|})^{|R|} = 2^{|R|^{2}} = 2^{|R|} = \gamma ,
\]
as required.  
\end{proof}

As a simple case of this lemma, we have the following corollary.  

\begin{corollary}\label{lem-size-bound}
Suppose $\cat{A}$ is a Grothendieck category, $X\in \Ch (\cat{A})$, and
$x\in X_{n}$ for some $n$.  Then the smallest subcomplex $Y$ of $X$ in
$\Ch (\cat{A})$ containing $x$ has $|Y|\leq \gamma $.
\end{corollary}

\begin{proof}
The smallest subcomplex $Y'$ of $R$-modules containing $x$ is simply
$R/\ann (x)$ in degree $n$ and $R/\ann (dx)$ in degree $n-1$, so
$|Y'|\leq \gamma $.  Since localization is exact, the localization $Y$
of $Y'$ will be a subcomplex in $\Ch (\cat{A})$ containing $x$.
Lemma~\ref{lem-cardinality} guarantees that $|Y|\leq \gamma $.
\end{proof}

We also need the following standard lemma.  Given a class of maps $K$,
$K\proj $ is the class of maps which look projective to $K$; that is,
they have the \llp $K$.  See~\cite[Section~2.1]{hovey-model} for the
precise definition.

\begin{lemma}\label{lem-K-proj}
Let $\cat{A}$ denote an abelian category with enough injectives.  Let
$K$ denote the class of surjections in $\Ch (\cat{A})$ whose kernel is
an injective object of $\Ch (\cat{A})$.  Then $K\proj $ is the class of
injections.  Furthermore, $(K\proj )\inj =K$.
\end{lemma}

\begin{proof}
We first show that any map in $K\proj $ is injective.  Recall the disk
functor $D^{n}\mathcolon \cat{A}\xrightarrow{}\Ch (\cat{A})$ that takes
an object $X$ to the complex which is $X$ in degrees $n$ and $n-1$, and
$0$ elsewhere.  The functor $D^{n}$ is right adjoint to the exact
functor $X\mapsto X_{n-1}$.  Thus $D^{n}(X)$ is injective whenever $X$
is injective in $\cat{A}$.  In particular, suppose $i\mathcolon
A\xrightarrow{}B$ is a map of complexes with kernel $C$.  Fix $n$, and
embed $C_{n}$ into an injective object $M$.  This embedding extends to a
map $A_{n}\xrightarrow{f}M$, and so defines a map of complexes
$A\xrightarrow{}D^{n+1}M$, which is $f$ in degree $n$.  This map
obviously cannot extend to a map $B\xrightarrow{}D^{n+1}M$ unless
$C_{n}$ is $0$.  Since the map $D^{n+1}M\xrightarrow{}0$ is in $K$, this
shows that every map in $K\proj $ is an injection.

Conversely, suppose we have a commutative diagram in $\Ch (\cat{A})$ as
follows,
\[
\begin{CD}
A @>>f> X \\
@ViVV @VVpV \\
B @>>g> Y
\end{CD}
\]
where $i$ is an injection and $p$ is a surjection with injective kernel
$W$.  Since $W$ is injective in $\Ch (\cat{A})$, there is a splitting
$q\mathcolon Y\xrightarrow{}X$ of $p$.  We have $p(qgi-f)=0$, so, since
$W$ is injective, there is an extension $h\mathcolon B\xrightarrow{}W$
such that $hi=qgi-f$.  Then $qg-h\mathcolon B\xrightarrow{}X$ is the
desired lift.  Hence $i$ is in $K\proj $.  

Now, we always have $(K\proj )\inj \supseteq K$.  Conversely, suppose
$p\mathcolon X\xrightarrow{}Y$ has the \rlp all injections.  Consider
the map $D^{n}(Y_{n})\xrightarrow{}Y$ that is the identity in degree
$n$.  Since $p$ has the \rlp all injections, there is a lift
$D^{n}Y\xrightarrow{}X$ of this map.  This shows that $p$ is a split
surjection in each dimension.  Since the map $\ker p\xrightarrow{}0$ is
a pullback of $p$, it too will have the \rlp all injections, and so
$\ker p$ is injective as an object of $\Ch (\cat{A})$.
\end{proof}

We can now prove that $I$ generates all injections. 

\begin{proposition}\label{prop-inj-cof}
Suppose $\cat{A}$ is a Grothendieck category.  The class $I\cof $ is the
class of injections, and the class $I\inj $ is the class of surjections
whose kernel is an injective object of $\Ch (\cat{A})$.
\end{proposition}

\begin{proof}
Let $K$ denote the class of surjections whose kernel is injective.
Applying Lemma~\ref{lem-K-proj}, we see that $I\subseteq K\proj $, so
$I\cof \subseteq (K\proj )\cof =K\proj $.  Thus $I\cof $ consists of
injections.  Furthermore, if we can show $I\cof $ consists of all
injections, Lemma~\ref{lem-K-proj} will show that $I\inj $ is $K$, as
required.

So suppose $i\mathcolon A\xrightarrow{}B$ is an injection.  To show that
$i\in I\cof $, we will show that $i$ has the \llp $I\inj $.  So
suppose $p\mathcolon X\xrightarrow{}Y$ is in $I\inj $, and we have a
commutative diagram as follows. 
\[
\begin{CD}
A @>f>> X \\
@ViVV @VVpV \\
B @>>g> Y
\end{CD}
\]
Let $T$ be the set of partial lifts of this diagram, so $T$ is the set
of all pairs $(C,h)$, where $C$ is a subcomplex of $B$ containing $i(A)$
and $h\mathcolon C\xrightarrow{}X$ is a chain map such that $ph=g|_{C}$
and $hi=f$.  Then $T$ is a partially ordered set, where $(C,h)\leq
(C',h')$ if $C'$ contains $C$ and $h'$ extends $h$.  The set $T$ is
nonempty and we claim that every chain in $T$ has an upper bound.
Indeed, given a chain $(C_{i},h_{i})$ in $T$, the colimit $C$ of the
$C_{i}$ is still a subcomplex of $B$, by the AB5 condition, and the
union of the $h_{i}$ defines a lift on $C$.  Thus there is a maximal
element $(M,h)$ of $T$.  Suppose that $M$ is not all of $B$, and choose
a homogeneous element $x\in B$ that is not in $M$.  Let $Z$ be the
smallest subcomplex (in $\Ch (\cat{A})$) of $B$ containing $x$, so that
$|Z|\leq \gamma $ by Corollary~\ref{lem-size-bound}.  Let $M'$ denote
the subcomplex of $B$ generated by $M$ and $x$, so that we have the
pushout diagram below.
\[
\begin{CD}
M\cap Z @>>> Z \\
@VVV @VVV \\
M @>>> M'
\end{CD}
\]
Since the top horizontal map is in $I$, the bottom horizontal map is in
$I\cof $.  Hence there is a lift $h'$ in the following diagram.  
\[
\begin{CD}
M @>h>> X \\
@VVV @VVpV \\
M' @>>g> Y
\end{CD}
\]
This lift violates the maximality of $(M,h)$, so we must have $M=B$.
Hence $i\in I\cof $, as required.  
\end{proof}

\begin{corollary}\label{cor-inj-cof}
Suppose $\cat{A}$ is a Grothendieck category.  Then every injective
object of $\Ch (\cat{A})$ is injectively fibrant and has no homology.
Every map in $I\inj $ is an injective fibration and a quasi-isomorphism.
\end{corollary}

\begin{proof}
The second statement follows from the first.  Indeed, a map in $I\inj $
has the \rlp all injections, so in particular is an injective fibration.
If $p\in I\inj $, then $\ker p$ is an injective object of $\Ch
(\cat{A})$, so has no homology by the first statement.  Thus $p$ is a
homology isomorphism, by the long exact sequence.  

Now suppose $X$ is an injective object of $\Ch (\cat{A})$.  Certainly
$X$ is injectively fibrant.  To see that $X$ has no homology, let $Y$
denote the complex defined by $Y_{n}=X_{n}\oplus X_{n-1}$ with
$d(x,y)=(dx+y,-dy)$.  Then $X\xrightarrow{}Y$ is an inclusion of
complexes, so since $X$ is injective, has a retraction
$Y\xrightarrow{}X$.  This retraction is equivalent to a contracting
homotopy of $X$, so in particular $X$ has no homology.
\end{proof}

To complete the proof of Theorem~\ref{thm-inj}, we must show that $J\cof
$ is the class of injective quasi-isomorphisms, from which it will
follow that $J\inj $ is the class of injective fibrations.  We begin
with the following crucial, but technical, lemma.

\begin{lemma}\label{lem-inj-weak-equivs-gamma}
Suppose $\cat{A}$ is a Grothendieck category.  Suppose $i\mathcolon
A\xrightarrow{}B$ is an injective quasi-isomorphism in $\Ch (\cat{A})$.
For every subcomplex $C$ of $B$ in $\Ch (\cat{A})$ with $|C|\leq \gamma
$, there is a subcomplex $D$ of $B$ in $\Ch (\cat{A})$ containing $C$
such that $|D|\leq \gamma $ and $i\mathcolon D\cap A\xrightarrow{}D$ is
a weak equivalence.
\end{lemma}

\begin{proof}
The failure of $C\cap A\xrightarrow{}C$ to be a quasi-isomorphism is
measured by $H_{*}(C/C\cap A)$.  Suppose for the moment that for every
homogeneous element $x$ of $H_{*}(C/C\cap A)$, we can find a subcomplex
$C(x)$ containing $C$ such that $|C(x)|\leq \gamma $ and the map
$H_{*}(C/C\cap A)\xrightarrow{}H_{*}(C(x)/C(x)\cap A)$ sends $x$ to $0$.
Since $|C|\leq \gamma $, the $R$-module homology of $C/C\cap A$ has size
$\leq \gamma $.  But then Lemma~\ref{lem-cardinality} assures that
$|H_{*}(C/C\cap A)|\leq \gamma $, so have $\leq \gamma $ choices for
$x$.  We can therefore take the union of all the $C(x)$ to form a new
subcomplex $FC$ with $|FC|\leq \gamma $ (using
Lemma~\ref{lem-cardinality} again), such that the induced map
$H_{*}(C/C\cap A)\xrightarrow{}H_{*}(FC/FC\cap A)$ is the zero map.

Now iterate this construction to form a sequence $F^{n}C$, and let $D$
be the colimit of all the $F^{n}C$.  Then $|D|\leq \gamma $, by
Lemma~\ref{lem-cardinality}.  Note that $D/D\cap A$ is the colimit of
the $F^{n}C/F^{n}C\cap A$, by commuting colimits.
Lemma~\ref{lem-colimits-homology} then shows that $H_{*}(D/D\cap A)=0$,
as required.

To complete the proof, we must construct the complex $C(x)$.  The
construction we give is fairly complicated; we do not know if there is a
simpler one.  Let us denote the $R$-module homology of a complex $X$ by
$\widetilde{H}(X)$ and let us denote the torsion submodule of an
$R$-module $M$ by $tM$.  Then the class $x$ is represented by a homomorphism
$f\mathcolon \ideal{a}\xrightarrow{}
\widetilde{H}_{n}(C/C\cap A)/t\widetilde{H}_{n}(C/C\cap A)$, for some
left ideal $\ideal{a}$ of $R$ such that $R/\ideal{a}$ is in the torsion
theory.  The class $x$ must map to $0$ in $H_{n}(B/A)$, since
$A\xrightarrow{}B$ is a quasi-isomorphism.  This means that there is a
subideal $\ideal{b}$ of $\ideal{a}$ with $R/\ideal{b}$ also in the
torsion theory, such that the composite
\[
\ideal{b} \xrightarrow{f} \widetilde{H}_{n}(C/C\cap
A)/t\widetilde{H}_{n}(C/C\cap A) \xrightarrow{}
\widetilde{H}_{n}(B/A)/t\widetilde{H}_{n}(B/A)
\]
is the zero map.  We need to construct $C(x)$ so that this map is
already the zero map when $C(x)$ replaces $B$. 

For each $y\in \ideal{b}$, choose an element $z_{y}$ in $C_{n}$ with
$dz_{y}\in A_{n}$ whose homology class $[z_{y}]$ is a representative for
$f(y)$.  Since the above composite is $0$, there is an ideal
$\ideal{c}_{y}$ such that $R/\ideal{c}_{y}$ is in the torsion theory,
and $\ideal{c}_{y}[z_{y}]=0$ in $\widetilde{H}_{n}(B/A)$.  This means
that, for every $w\in \ideal{c}_{y}$, there is an element $v_{w,y}\in
B_{n+1}$ such that $wz_{y}-dv_{w,y}\in A_{n}$.  We define $C(x)$ to be
the smallest subcomplex of $B$ containing $C$ and all the $v_{w,y}$.  It
is clear from the construction that $x$ goes to $0$ in
$H_{n}(C(x)/C(x)\cap A)$.  Since there are $\leq |R|\leq \gamma $
choices for $w$ and $y$, the smallest subcomplex of $R$-modules
containing $C$ and the $v_{w,y}$ has size $\leq \gamma $.
Lemma~\ref{lem-cardinality} then shows $|C(x)|\leq \gamma $, as
required.
\end{proof}

With this lemma in hand, it is now not difficult to show that $J\cof $
is the class of injective quasi-isomorphisms.  

\begin{proposition}\label{prop-inj-triv-cof}
Suppose $\cat{A}$ is a Grothendieck category.  Then the class $J\cof $
consists of the injective quasi-isomorphisms, and the class $J\inj $
consists of the injective fibrations.
\end{proposition}

By the recognition theorem~\cite[Theorem~2.1.19]{hovey-model}, this
proposition completes the proof of Theorem~\ref{thm-inj}.  

\begin{proof}
The second statement is an immediate corollary of the first.  By
Corollary~\ref{cor-homology-J}, the maps of $J\cof $ are injective
quasi-isomorphisms.  Now suppose $i\mathcolon A\xrightarrow{}B$ is an
injective quasi-isomorphism.  To show that $i\in J\cof $, we show that
$i$ has the \llp $J\inj $.  So suppose $p$ is in $J\inj $, and we have a
commutative diagram as follows.
\[
\begin{CD}
A @>f>> X \\
@ViVV @VVpV \\
B @>>g> Y
\end{CD}
\]
Let $T$ denote the set of partial lifts $(C,h)$, where $C$ is a
subcomplex of $B$ containing $iA$ such that the map $i\mathcolon
A\xrightarrow{}C$ is a quasi-isomorphism, and $h\mathcolon
C\xrightarrow{}X$ is a partial lift in our diagram.  Then $T$ is
obviously partially ordered and nonempty.
Proposition~\ref{prop-colimits-weak-equivs} and the argument used in the
proof of Proposition~\ref{prop-inj-cof} imply that a chain in $T$ has an
upper bound.  Zorn's lemma then gives us a maximal element $(M,h)$ of
$T$.  Suppose $M$ is not all of $B$, and choose an element $x$ in $B$
but not in $M$.  Let $C$ denote the subcomplex of $B$ generated by $x$,
so $|C|\leq \gamma $ by Corollary~\ref{lem-size-bound}.  Since
$M\xrightarrow{}B$ is a quasi-isomorphism,
Lemma~\ref{lem-inj-weak-equivs-gamma} implies that there is a complex
$D$ containing $C$ such that $|D|\leq \gamma $ and the map $D\cap
M\xrightarrow{}D$ is a quasi-isomorphism.  Let $N$ denote the subcomplex
of $B$ generated by $M$ and $D$.  Then the map $M\xrightarrow{}N$ is in
$J\cof $, since it is a pushout of $D\cap M\xrightarrow{}D$.  Since
$p\in J\inj $, there is an extension of $h$ to $N$, contradicting the
maximality of $(M,h)$.  Therefore we must have had $M=B$, and so $i\in
J\cof $, as required.
\end{proof}

To complete the description of the injective model structure, we would
like to characterize the injective fibrations.  This characterization is
precisely the same as the corresponding characterization in the category
of chain complexes of modules, found
in~\cite[2.3.16--20]{hovey-model}, with the same proofs. 

\begin{proposition}\label{prop-inj-fibs}
Suppose $\cat{A}$ is a Grothendieck category.  Then a map $p\in \Ch
(\cat{A})$ is an injective fibration if and only if it is a split
surjection in each degree with injectively fibrant kernel.  Any
injectively fibrant complex is a complex of injective objects, and any
bounded above complex of injective objects is injectively fibrant.
\end{proposition}

We now discuss the functoriality of the injective model structure.

\begin{proposition}\label{prop-inj-functor}
Suppose $F\mathcolon \cat{A}\xrightarrow{}\cat{B}$ is a functor between
Grothendieck categories, with right adjoint $U$.  Then $F$ induces a
Quillen adjunction $F\mathcolon \Ch (\cat{A})\xrightarrow{}\Ch
(\cat{B})$ between the injective model structures if and only if $F$ is
exact.  
\end{proposition}

\begin{proof}
If $F$ is exact, then clearly $F$ preserves injections and all
quasi-isomorphisms, so preserves cofibrations and trivial cofibrations.
Conversely, suppose $F$ preserves cofibrations and trivial
cofibrations.  We can think of an exact sequence as a complex $X$ with
no homology, so $0\xrightarrow{}X$ is a trivial cofibration.  Then
$0\xrightarrow{}FX$ must also be a trivial cofibration, so $F$ must be
exact.  
\end{proof}

Note that $F$ being exact is equivalent to $U$ being additive and
preserving injectives.  This proposition is expected, but not very
satisfying.  It means we cannot use the injective model structure to
form any interesting total left derived functors, since such a total
left derived functor would be defined by first replacing an object by a
cofibrant object weakly equivalent to it, and every object is already
cofibrant.  We can use the injective model structure to form some right
derived functors.

In fact, we can use it to form more right derived functors than one
might first expect.  The construction of the total right derived functor
of $U$~\cite[Definition~1.3.6]{hovey-model} only requires that $U$
preserve weak equivalences between injectively fibrant objects.  But a
quasi-isomorphism between injectively fibrant chain complexes is in fact
a chain homotopy equivalence (since every object is cofibrant).  So in
order to construct the right derived functor of $U$, we only need to
insure that $U$ preserve chain homotopy.  For this, all we require is
that $U$ be additive.  

We have then proved the following proposition. 

\begin{proposition}\label{prop-inj-right}
Suppose $U\mathcolon \cat{A}\xrightarrow{}\cat{B}$ is an additive
functor between Grothendieck categories.  Then the total right derived
functor $RU\mathcolon \cat{D}(\cat{A})\xrightarrow{}\cat{D}(\cat{B})$ of
$U$ exists.  
\end{proposition}

This recovers the usual right derived functors of $U$: if $X\in
\cat{A}$, we have $(R^{i}U)(X)=H_{i}((RU)X)$.  The functor $RU$ is of
course calculated by replacing $X$ by an injective resolution (or an
injectively fibrant approximation if $X$ is a complex), then applying
$U$.  In particular, if $f$ is a map of ringed spaces, we recover the
total right derived functor $Rf_{*}$ between complexes of sheaves in
this way.

\section{An alternative approach}\label{sec-alt}

We have already discussed the drawbacks of the injective model structure
on a Grothendieck category $\cat{A}$.  In this section, we offer another
approach; we will apply it to the category of sheaves on a ringed space
satisfying a hypthesis related to finite global dimension in the next
section.  Though this is our only application of this approach, we
present the method in a general fashion in the hope that it may find
other applications.  This approach is based on the standard projective
model structure when $\cat{A}=R\Mod $ for some ring $R$, and
generalizations of it considered by Christensen
in~\cite{christensen-derived}.  Recall
from~\cite[Section~2.3]{hovey-model} that the projective model structure
on $\Ch (\cat{A})$, where $\cat{A}=R\Mod $ for some ring $R$, is a
cofibrantly generated model structure, with generating cofibrations
$I=\{S^{n-1}R\xrightarrow{}D^{n}R \}$ and generating trivial
cofibrations $J=\{0\xrightarrow{}D^{n}R \}$.  Here $n$ runs through all
integers, $S^{n-1}M$ is the complex whose only nonzero object is $M$ in
dimension $n-1$, and $D^{n}M$ is the complex whose only nonzero objects
are $M$ in dimensions $n$ and $n-1$.  Our plan is to replace the map
$0\xrightarrow{}R$ by a set of monomorphisms $\cat{M}$.



\begin{definition}\label{defn-flat-I}
Suppose $\cat{M}$ is a set of monomorphisms in a Grothendieck category
$\cat{A}$.  Let $\cat{F}$ denote the set of codomains of the maps of
$\cat{M}$.  We will say that $\cat{M}$ is \emph{pointed} if $0\in
\cat{F}$ and, if $F\in \cat{F}$, then $0\xrightarrow{}F$ is in
$\cat{M}$.  Define $J$ to be the set of all $D^{n}f$, where $n$ is an
integer and $f\in \cat{M}$.  Define $I$ to be the union of $J$ and the
maps $S^{n-1}F\xrightarrow{}D^{n}F$ for $F\in \cat{F}$ and $n$ an
integer.  Then define a map $p$ to be a \emph{$\cat{M}$-fibration} if
$p$ is in $J\inj $, define $p$ to be a \emph{$\cat{M}$-cofibration} if
$p$ is in $I\cof $. 
\end{definition}

If $\cat{M}$ consists only of the maps $0\xrightarrow{}F$
for $F\in \cat{F}$, then we recover the definitions
of~\cite{christensen-derived}.


Our goal is to determine conditions on $\cat{M}$ under which the
quasi-isomorphisms, the $\cat{M}$-cofibrations, and the
$\cat{M}$-fibrations determine a model structure on $\Ch (\cat{A})$.  We
use the recognition theorem~\cite[Theorem~2.1.19]{hovey-model}.  Since
the maps of $J$ are injective quasi-isomorphisms in $I\cof $, the maps
of $J\cof $ will also be, by Corollary~\ref{cor-homology-J}.  Hence we
need to show that the maps of $I\inj $ coincide with the maps that are
both $\cat{M}$-fibrations and quasi-isomorphisms.  





We begin by characterizing the $\cat{M}$-fibrations.  

\begin{definition}\label{defn-flasque}
Suppose $\cat{M}$ is a pointed set of monomorphisms in a Grothendieck
category $\cat{A}$.  Define an object $X$ of $\cat{A}$ to be
\emph{$\cat{M}$-flasque} if $\cat{A}(f,X)$ is surjective for all $f\in
\cat{M}$.
\end{definition}

This definition is a generalization of the usual notion of flasque, or
flabby, sheaves.  We will discuss this in detail in the next section. 

Let us denote the category of chain complexes of abelian groups by $\Ch
 (\Z )$.  

\begin{proposition}\label{prop-fibs}
Suppose $\cat{M}$ is a pointed set of monomorphisms in a Grothendieck
category $\cat{A}$.  Then a map $p\mathcolon X\xrightarrow{}Y$ in $\Ch
(\cat{A})$ is a $\cat{M}$-fibration if and only if $\cat{A}(F,p)$ is a
surjection in $\Ch (\Z )$ for all $F$ in $\cat{F}$ and $\ker p$ is
dimensionwise $\cat{M}$-flasque.
\end{proposition}

In particular, if $\cat{F}$ is a set of generators for $\cat{A}$, then
$\cat{M}$-fibrations are surjective.  To see this, consider the map from
$Y_{n}$ into the cokernel of $p_{n}$.  

\begin{proof}
Adjointness implies that $p$ has the \rlp
$D^{n}B\xrightarrow{D^{n}f}D^{n}C$ if and only if the map
\[
\cat{A}(C,X_{n})\xrightarrow{}\cat{A}(C,Y_{n})\times _{\cat{A}(B,Y_{n})}
\cat{A}(B,X_{n}) 
\]
is surjective.  Applying this when $f$ is the map $0\xrightarrow{}F$ for
$F\in \cat{F}$, we find that, if $p$ is a $\cat{M}$-fibration, then
$\cat{A}(F,p)$ is surjective.  Furthermore, if $p$ is a
$\cat{M}$-fibration, then $\ker p\xrightarrow{}0$ is in $J\inj $.  
Applying the above criterion, we find that $\ker p$ is dimensionwise
$\cat{M}$-flasque.  

Conversely, suppose $\cat{A}(F,p)$ is a surjection for all $F\in
\cat{F}$ and $K=\ker p$ is dimensionwise $\cat{M}$-flasque.  Suppose
$f\mathcolon B\xrightarrow{}C$ is in $\cat{M}$.  We have an exact
sequence
\[
0\xrightarrow{}\cat{A}(B,K_{n}) \xrightarrow{} \cat{A}(B,X_{n})
\xrightarrow{}\cat{A}(B,Y_{n})
\]
and a similar exact sequence that is in fact short exact when $B$ is
replaced by $C$.  By pulling back the exact sequence for $B$ through
the map $\cat{A}(f,Y_{n})$, we obtain the following commutative diagram
whose top row is short exact and whose bottom row only misses being
short exact because the right map is not necessarily surjective.  
\[
\begin{CD}
\cat{A}(C,K_{n}) @>>> \cat{A}(C,X_{n}) @>>> \cat{A}(C,Y_{n}) \\
@V\cat{A}(f,K_{n})VV @VVV @| \\
\cat{A}(B,K_{n}) @>>> \cat{A}(C,Y_{n})\times _{\cat{A}(B,Y_{n})}
\cat{A}(B,X_{n}) @>>> \cat{A}(C,Y_{n})
\end{CD}
\]
Since $K$ is dimensionwise $\cat{M}$-flasque, the left-hand vertical map
is surjective.  A standard diagram chase, as in the snake lemma, then
show that the middle vertical map is surjective, so $p$ is a
$\cat{M}$-fibration.
\end{proof}

\begin{proposition}\label{prop-trivial-fibs}
Suppose $\cat{M}$ is a pointed set of monomorphisms in a Grothendieck
category $\cat{A}$.  Suppose in addition that the set $\cat{F}$ of
codomains of $\cat{M}$ generates $\cat{A}$.  Then every map of complexes
$p\mathcolon X\xrightarrow{}Y$ in $I\inj $ is both a $\cat{M}$-fibration
and a quasi-isomorphism.
\end{proposition}

\begin{proof}
Recall that the functor $S^{n-1}\mathcolon \cat{A}\xrightarrow{}\Ch
(\cat{A})$ is left adjoint to the functor that takes $X$ to $Z_{n-1}X$,
the cycles in $X_{n-1}$.  This implies that $p$ is in $I\inj $ if and
only if it is a $\cat{M}$-fibration and the map
\[
\cat{A}(F,X_{n})\xrightarrow{}\cat{A}(F,Y_{n})\times_{\cat{A}(F,Z_{n-1}Y)}
\cat{A}(F,Z_{n-1}X)
\]
is surjective for all $n$ and $F\in \cat{F}$.
Let $K=\ker p$.  If $p\in I\inj $, then the map $K\xrightarrow{}0$ is as
well.  Hence the map $\cat{A}(F,K_{n})\xrightarrow{}\cat{A}(F,Z_{n-1}K)$ is
surjective for all $n$ and all $F\in \cat{F}$.  Since $\cat{F}$ is a set
of generators for $\cat{A}$, this implies that the map
$K_{n}\xrightarrow{}Z_{n-1}K$ is surjective, and hence that $K$ has no
homology.  A similar argument shows that $p$ is surjective, and so the
long exact sequence implies that $p$ is a quasi-isomorphism.    
\end{proof}

If $\cat{F}$ is not a generating set for $\cat{A}$, we can still say
that, if $p\in I\inj $, then $\cat{A}(F,p)$ is a surjective
quasi-isomorphism for all $F\in \cat{F}$.

To complete the construction of our model structure, we need to know
that every map that is both a $\cat{M}$-fibration and a
quasi-isomorphism is in $I\inj $.  We begin with a lemma.

\begin{lemma}\label{lem-kernel}
Suppose $\cat{M}$ is a pointed set of monomorphisms in a Grothendieck
category $\cat{A}$, and let $\cat{F}$ be the set of codomains of
$\cat{M}$.  Suppose $p\mathcolon X\xrightarrow{}Y$ is a map in $\Ch
(\cat{A})$ such that $\cat{A}(F,p)$ is surjective for all $F\in
\cat{F}$.  Then $p$ is in $I\inj $ if and only if $\ker
p\xrightarrow{}0$ is in $I\inj $.
\end{lemma}

\begin{proof}
The only if implication is clear.  Suppose $\cat{A}(F,p)$ is surjective
for all $F\in \cat{F}$, and let $K=\ker p$.  Suppose $K\xrightarrow{}0$
is in $I\inj $.  In particular, this means that $K$ is dimensionwise
$\cat{M}$-flasque, so $p\in J\inj $.  In order to show that $p$ is in
$I\inj $, we must show that, given $F\in \cat{F}$, a map $x\mathcolon
F\xrightarrow{}Z_{n-1}X$, and a map $y\mathcolon F\xrightarrow{}Y_{n}$
such that $d\circ y=p\circ x$, there is a map $x'\mathcolon
F\xrightarrow{}X_{n}$ such that $p\circ x'=y$ and $d\circ x'=x$.  First
choose $z\mathcolon F\xrightarrow{}X_{n}$ such that $p\circ z=y$, using
the fact that $\cat{A}(F,p)$ is surjective.  Then $p\circ (d\circ
z-x)=0$, so $dz-x \mathcolon F\xrightarrow{} Z_{n-1}K$.  Since
$K\xrightarrow{}0$ is in $I\inj $, there is a map $w\mathcolon
F\xrightarrow{}K_{n}$ such that $d\circ w=d\circ z-x$.  Now let
$x'=z-w$.
\end{proof}

\begin{proposition}\label{prop-flat-triv-fib}
Suppose $\cat{M}$ is a pointed set of monomorphisms in a Grothendieck
category $\cat{A}$.  Suppose that the set $\cat{F}$ of codomains of
$\cat{M}$ generates $\cat{A}$, and, furthermore, suppose that if $K$ is
an acyclic, dimensionwise $\cat{M}$-flasque, complex and $F\in \cat{F}$,
then $\cat{A}(F,K)$ is an acyclic complex of abelian groups.  Then, if
$p\mathcolon X\xrightarrow{}Y$ is a $\cat{M}$-fibration and
quasi-isomorphism in $\Ch (\cat{A})$, then $p$ is in $I\inj $.
\end{proposition}

\begin{proof}
By Lemma~\ref{lem-kernel}, it suffices to show that $K=\ker
p\xrightarrow{}0$ is in $I\inj $.  But $K$ is an acyclic
dimensionwise flasque complex, and so $\cat{A}(F,K)$ is also acyclic.
Hence the map $\cat{A}(F,K_{n})\xrightarrow{}\cat{A}(F,Z_{n-1}K)$ is
surjective, and so $K\xrightarrow{}0$ is in $I\inj $.  
\end{proof}

We have proved the following theorem.  

\begin{theorem}\label{thm-alt-model}
Suppose $\cat{M}$ is a pointed set of monomorphisms in a Grothendieck
category $\cat{A}$ such that the set of codomains $\cat{F}$ of $\cat{M}$
forms a generating set of $\cat{A}$ and, for all acyclic, dimensionwise
$\cat{M}$-flasque, complexes $X$ and for all $F\in \cat{F}$, the complex
$\cat{A}(F,X)$ is acyclic.  Then $\Ch (\cat{A})$ is a proper cofibrantly
generated model category, where the weak equivalences are the
quasi-isomorphisms, the fibrations are the $\cat{M}$-fibrations, and the
cofibrations are the $\cat{M}$-cofibrations.
\end{theorem}

One interesting feature of the hypotheses of this theorem is that, if
they are true for a given set of monomorphisms $\cat{M}$ with codomains
$\cat{F}$, then they remain true if we expand $\cat{M}$ by adding any
set of monomorphisms whose codomains are all in $\cat{F}$.  So in fact
we get many different model structures with the same weak equivalences,
all relying on more or less stringent definitions of ``flasque''.  

One might hope that we would still get a model structure on $\Ch
(\cat{A})$ if we drop all hypotheses about the set of monomorphisms
$\cat{M}$.  The weak equivalences would have to change, probably to maps
$f$ such that $\cat{A}(F,f)$ is a quasi-isomorphism for all $F\in
\cat{F}$.  With this definition, an appropriately modified version of
Proposition~\ref{prop-flat-triv-fib} does hold.  However, we do not know
if the maps of $J\cof $ are weak equivalences with this definition.  

\section{Generators of finite projective dimension}\label{sec-fin-dim}

In this section, we apply the method of the previous section to
construct a new model structure on $\Ch (\cat{A})$, when $\cat{A}$ is a
Grothendieck category with generators of finite projective dimension.  
Recall that an object $B$ is said to have \emph{finite projective
dimension} if there is an integer $n_{0}$ such that $\Ext
_{\cat{A}}^{n}(B,C)=0$ for all $n\geq n_{0}$ and all object $C$ of
$\cat{A}$.  One normally thinks of an object of finite projective
dimension as being the $0$th homology group of a finite complex of
projectives, but this will not be true unless there are enough
projectives in the category.  In the categories we are interested in,
this is almost never true.  

Nevertheless, objects of finite projective dimension are useful in
constructing a model structure because of the following lemma.  

\begin{lemma}\label{lem-fin-dim}
Suppose $\cat{A}$ is a Grothendieck category, $F\in \cat{A}$ has finite
projective dimension, and $X\in \Ch (\cat{A})$ is an acyclic complex
such that $\Ext _{\cat{A}}^{i}(F,X_{n})=0$ for all $i>0$ and all $n$.
Then $\cat{A}(F,X)$ is still acyclic.  
\end{lemma}

\begin{proof}
Since $X$ is acyclic, we have a short exact sequence 
\[
Z_{n}X \xrightarrow{} X_{n} \xrightarrow{}Z_{n-1}X .
\]
Since $\Ext _{\cat{A}}^{i}(F,X_{n})=0$ for all $i>0$, this gives us an
exact sequence 
\[
0 \xrightarrow{} \cat{A}(F,Z_{n}X) \xrightarrow{} \cat{A}(F,X_{n})
\xrightarrow{} \cat{A}(F,Z_{n-1}X) \xrightarrow{} \Ext
^{1}_{\cat{A}}(F,Z_{n}X) \xrightarrow{} 0
\]
and isomorphisms $\Ext ^{i}_{\cat{A}}(F,Z_{n-1}X)\cong \Ext
^{i+1}_{\cat{A}}(F,Z_{n}X)$ for $i>0$.  Thus 
\[
\Ext ^{1}_{\cat{A}}(F,Z_{n}X) \cong \Ext ^{m+1}_{\cat{A}}(F,Z_{m+n}X)
\]
for all $m\geq 0$.  Since $F$ has finite projective dimension, this
implies $\Ext ^{1}_{\cat{A}}(F,Z_{n}X)=0$ for all $n$.  It follows
that $\cat{A}(F,X)$ is acyclic.  
\end{proof}

\begin{theorem}\label{thm-fin-dim}
Suppose $\cat{A}$ is a Grothendieck category with a set of generators
$\cat{F}$, each element of which has finite projective dimension.  Let
$\cat{M}$ denote the set of inclusions $A\xrightarrow{}F$ of subobjects
of objects $F\in \cat{F}$.  Then there is a proper cofibrantly generated
model structure on $\Ch (\cat{A})$, where the weak equivalences are the
quasi-isomorphisms, the fibrations are the dimensionwise split
surjections with dimensionwise injective kernel, and the cofibrations
are the $\cat{M}$-fibrations.  
\end{theorem}

\begin{proof}
Note first that the $\cat{M}$-flasque objects of $\cat{A}$ coincide with
the injective objects, by~\cite[Prop.~V.2.9]{stenstrom}.
Lemma~\ref{lem-fin-dim} implies that if $X$ is an acyclic, dimensionwise
injective, complex, then $\cat{A}(F,X)$ is acyclic for all $F\in
\cat{F}$.  Hence Theorem~\ref{thm-alt-model} gives us a model
structure.  Any $\cat{M}$-fibration is a surjection with dimensionwise
injective kernel, by Proposition~\ref{prop-fibs}, and therefore must be
a dimensionwise split surjection.  Conversely, a dimensionwise split
surjection with dimensionwise injective kernel certainly satisfies the
conditions of Proposition~\ref{prop-fibs}, so is a $\cat{M}$-fibration.  
\end{proof}

This model structure is related to the injective model structure; the
identity functor from this model structure to the injective model
structure is a Quillen equivalence.  It appears to be new even when
$\cat{A}$ is the category of modules over a ring $R$.  The generating
cofibrations and trivial cofibrations in this model structure are
explicit, and the fibrations are easier to understand than the injective
fibrations.  On the other hand, we know nothing about the cofibrations
in this model structure.  

In general, this model structure is poorly behaved with respect to
functors of abelian categories.  If $F$ is an additive functor with
right adjiont $U$, then $U$ will preserve fibrations in this model
structure if an only if $U$ preserves injectives, which is equivalent to
$F$ being exact.  But this is not enough to conclude that $F$ induces a
Quillen functor; we must also know that $U$ preserves acyclic complexes
of injectives.  This will happen if $U$ is exact, but may happen in some
other cases as well.  

We now consider an interesting example of this model structure.  Suppose
$S$ is a noetherian scheme.  We say that $S$ \emph{has enough locally
frees} if every coherent sheaf on $S$ is a quotient of a locally free
sheaf of finite rank.  For example, a noetherian, integral, separated,
locally factorial scheme has enough locally frees by a result of
Kleiman~\cite[Ex.~III.6.8]{hartshorne}.

\begin{proposition}\label{prop-enough}
Suppose $S$ is a noetherian scheme with enough locally frees.  In
addition, suppose that either $S$ is finite-dimensional or is separated.
Then the locally free sheaves of finite rank are generators of finite
projective dimension for the category $\QCo (S)$ of quasi-coherent
sheaves on $X$.
\end{proposition}

\begin{proof}
We first show that the locally frees generate $\QCo (S)$.
Deligne~\cite[Appendix, Prop.~2]{hartshorne-residues} shows that every
quasi-coherent sheaf is a colimit of finitely presented sheaves.  On a
noetherian scheme, finitely presented sheaves are coherent, and thus,
since $S$ has enough locally frees, are quotients of locally free
sheaves of finite rank.

Now let $F$ be a locally free sheaf of finite rank, and $C$ a
quasi-coherent sheaf of $\cat{O}$-modules on $S$.  By the corollary
to~\cite[Prop.~4.2.3]{grothendieck-tohoku}, we have
\[
\Ext ^{i}_{\cat{O}\Mod }(F,C) \cong H^{i}(S;\Hom (F,C))
\]
where $\Hom $ denotes sheaf Hom and the cohomology groups are sheaf
cohomology.  If $S$ is finite-dimensional, we can apply Grothendieck's
vanishing theorem~\cite[Theorem~III.2.7]{hartshorne} to conclude that
these cohomology groups are $0$ for large enough $i$.  If $S$ is
separated, then we can apply~\cite[Ex.~III.4.8]{hartshorne} to reach the
same conclusion, using the fact that $\Hom (F,C)$ is quasi-coherent.

This does not complete the proof, because these are $\Ext $ groups in
$\cat{O}\Mod $ rather than in $\QCo (S)$.  However, these two possibly
different $\Ext $ groups in fact coincide, because the exact inclusion
functor $\QCo (S)\xrightarrow{}\cat{O}\Mod $ has a right adjoint and
left inverse $Q$~\cite[p.~187]{SGA6} whenever $S$ is quasi-compact and
quasi-separated, as any noetherian scheme is.  In detail, given a
quasi-coherent sheaf $C$, we can first take an injective resolution
$I_{*}$ of $C$ in $\cat{O}\Mod $ and apply $Q$ to get a complex of
injectives $QI_{*}$ in $\QCo (S)$.  We claim that $QI_{*}$ is still
exact.  To see this, consider the short exact sequence
\[
0 \xrightarrow{} C \xrightarrow{} I_{0} \xrightarrow{} ZI_{1}
\xrightarrow{} 0 .
\]
After we apply $Q$, we get a long exact sequence involving the derived
functors $R^{i}Q$ of $Q$.  However, $R^{i}QC=0$ for $i>0$, by the last
paragraph of~\cite[p.~189]{SGA6}.  Furthermore, $R^{i}QI_{0}=0$ for
$i>0$ because $I_{0}$ is injective.  It follows that $(R^{i}Q)ZI_{1}=0$
for $i>0$ as well.  Repeating this argument on the short exact sequence
\[
0 \xrightarrow{} ZI_{1} \xrightarrow{} I_{1} \xrightarrow{} ZI_{2}
\xrightarrow{} 0,
\]
we find that $(R^{i}Q)ZI_{2}=0$ for $i>0$, and, by induction, that
$(R^{i}Q)ZI_{m}=0$ for all $m$ and $i>0$.  Hence $QI_{*}$ is still exact,
and so is an injective resolution of $C$ in $\QCo (S)$.  

Applying $\QCo (S)(B,-)$ to $QI_{*}$ and using adjointness, we find
that, if $B$ and $C$ are both quasi-coherent, then $\Ext ^{i}_{\QCo
(S)}(B,C)=\Ext ^{i}_{\cat{O}\Mod }(B,C)$, completing the proof.  
\end{proof}

Hence, as a corollary to Proposition~\ref{prop-enough} and
Theorem~\ref{thm-fin-dim}, we get the following theorem. 

\begin{theorem}\label{thm-quasi-coherent}
Suppose $S$ is a noetherian scheme with enough locally frees, and
suppose that $S$ is either finite-dimensional or separated.  Then there
is a proper, cofibrantly generated, model structure on the category
$\qch (S)$ of unbounded complexes of quasi-coherent sheaves, where the
weak equivalences are the quasi-isomorphisms and the fibrations are the
dimensionwise split surjections with dimensionwise injective kernel.
\end{theorem}

Let us call this model strucure the \emph{locally free model structure}.
We do not understand the cofibrations in the locally free model
structure, though we point out that $S^{n}F$ is cofibrant for any
locally free $F$, and $D^{n}A$ is cofibrant for any coherent sheaf $A$.
If $f\xrightarrow{}S\xrightarrow{}T$ is a map between schemes satisfying
the hypotheses of Theorem~\ref{thm-quasi-coherent}, then the functor
$f^{*}\mathcolon \QCo (T)\xrightarrow{}\QCo (S)$ will induce a Quillen
functor between the locally free model structures if and only if $f^{*}$
is exact; we have already seen that this is necessary, and it is
sufficient since $f^{*}$ preserves locally free sheaves of finite rank.

Despite these drawbacks, the locally free model structure does gives
some information about the derived category $D(\QCo (S))$.  

\begin{corollary}\label{cor-quasi-coherent}
Suppose $S$ is a noetherian scheme with enough locally frees, and either
$S$ is finite-dimensional or separated.  Then the locally free sheaves
of finite rank form a set of small weak generators for the derived
category $D(\QCo (S))$.  
\end{corollary}

\begin{proof}
The fact that the locally free sheaves form a set of weak generators
follows from~\cite[Section~7.3]{hovey-model}.  To see that they are
small, in the triangulated sense, we use the result
of~\cite[Section~7.4]{hovey-model}.  We must then show that, if $F$ is a
locally free sheaf of finite rank, the functor $\QCo (S)(F,-)$ preserves
all transfinite compositions.  Since we are on a noetherian scheme, we
can take the transfinite composition in the category of
presheaves~\cite[Ex.~II.1.11]{hartshorne}.  It is then easy to check the
desired result.  
\end{proof}

In case $S$ is a quasi-compact, quasi-separated scheme, we can use the
right adjoint $Q$ to the inclusion $\QCo (S)\xrightarrow{}\cat{O}\Mod $
to show that $\QCo (S)$ is a closed symmetric monoidal category under
the tensor product.  Thus $\qch (S)$ is also a closed symmetric monoidal
category.  It would be preferable, then, to have a model structure on
$\qch (S)$ that is compatible with the closed symmetric monoidal
structure, in the sense of~\cite[Chapter~4]{hovey-model}.  This
compatiblity condition is discussed before Theorem~\ref{thm-boxprod}
below.  Unfortunately, the locally free model structure is not
compatible with the tensor product.  

Despite this, it is known that $D(\QCo (S))$ is a symmetric monoidal
triangulated category, at least when $S$ is a finite-dimensional
noetherian scheme.  Indeed, Lipman~\cite[Section~2.5]{lipman}shows that
$D(\cat{O}\Mod )$ is a symmetric monoidal triangulated category.  But
$D(\QCo (S))$ is equivalent to the full subcategory of $D(\cat{O}\Mod )$
consisting of complexes with quasi-coherent cohomology, when $S$ is a
finite-dimensional noetherian scheme, by~\cite[p.~191]{SGA6}, and the
inclusion $D(\QCo (S))\xrightarrow{}D(\cat{O}\Mod )$ has a right adjoint
given by the right derived functor of $Q$.  It follows from this that
$D(\QCo (S))$ is a symmetric monoidal triangulated category.  

Furthermore, locally free sheaves of finite rank $F$ are strongly
dualizable in $D(\QCo (S))$.  Recall that this means that the natural
map 
\[
\Hom (F,\cat{O}) \otimes X \xrightarrow{}\Hom (F,X)
\]
is an isomorphism, where of course both the $\Hom $ and the tensor have
to be interpreted in $D(\QCo (S))$, so are really derived versions.
This follows from the corresponding fact in $\cat{O}\Mod $ itself, and
the fact that locally free sheaves are flat.  

In the language of~\cite{hovey-axiomatic}, then, we have proved the
following corollary.  

\begin{corollary}\label{cor-unital-alg}
Suppose $S$ is a finite-dimensional noetherian scheme with enough
locally frees.  Then the category $D(\QCo (S))$ is an unital algebraic
stable homotopy category, where the generators are the locally free
sheaves of finite rank.
\end{corollary}

\section{The flat model structure on sheaves}\label{sec-flat}

In this section, we apply the method of Theorem~\ref{thm-alt-model} to
the category $\cat{O}\Mod $ of sheaves over a ringed space
$(S,\cat{O})$.  In this case, there is a standard set of generators;
namely, the sheaves $\cat{O}_{U}$ for $U$ an open set of $S$.  Recall
that $\cat{O}_{U}$ is the sheafification of the presheaf that assigns
$V$ to $\cat{O}(V)$ if $V\subseteq U$, and to $0$ otherwise.  The stalk
of $\cat{O}_{U}$ at $x$ is $0$ if $x\not \in U$, and is $\cat{O}_{x}$ if
$x\in U$.  We have $\cat{O}\Mod (\cat{O}_{U},X)\cong X(U)$, which
implies easily that the $\cat{O}_{U}$ form a generating set for
$\cat{O}\Mod $.

Note that, if $V\subseteq U$, there is a natural monomorphism
$\cat{O}_{V}\xrightarrow{}\cat{O}_{U}$ corresponding to $1\in
\cat{O}_{U}(V)$.  Thus, we take the set of monomorphisms $\cat{M}$ of
the previous section to consist of these natural monomorphisms.  One can
then easily check that a sheaf $X$ is $\cat{M}$-flasque if and only if the
restriction maps $X(U)\xrightarrow{}X(V)$ are surjective whenever
$V\subseteq U$, corresponding to the usual notion of a flasque sheaf.  

To apply Theorem~\ref{thm-alt-model} we need to know that, if $X$ is an
acyclic complex of flasque sheaves, then $\cat{O}\Mod (\cat{O}_{U},X)$
is also acyclic; \ie that $X$ is acyclic as a complex of
\emph{presheaves}.  Unfortunately, this need not always be true.  Amnon
Neeman has constructed a complex $X$ of injective sheaves on
infinite-dimensional real projective space whose sheaf cohomology is
trivial, but whose presheaf cohomology is non-trivial.  The example is a
bit complicated, but is closely related to the example
in~\cite[Remark~2.3.18]{hovey-model}.

We therefore need a hypothesis on our ringed space to apply
Theorem~\ref{thm-alt-model}.  

\begin{definition}\label{defn-dim}
Define a ringed space $(S,\cat{O})$ to have \fgd \ if there is an integer
$n>0$ such that the sheaf cohomology $H^{n}(X)=0$ for all
$\cat{O}$-modules $X$.  Define $(S,\cat{O})$ to have \fhgd \  if every open
ringed subspace $(U,\cat{O}|_{U})$ has \fgd .
\end{definition}

We then get the following theorem.  

\begin{theorem}\label{thm-sheaves-model}
Suppose $(S,\cat{O})$ is a ringed space with \fhgd.  Then there is a
cofibrantly generated proper model structure on $\Ch (\cat{O}\Mod )$,
called the \emph{flat model structure}, where the weak equivalences are
the quasi-iso\-morph\-isms and the fibrations are the surjections with
dimensionwise flasque kernel.  
\end{theorem}

\begin{proof}
We apply Theorem~\ref{thm-alt-model}, taking the set $\cat{M}$ to be the
canonical inclusions $\cat{O}_{V}\xrightarrow{}\cat{O}_{U}$.  We use
Lemma~\ref{lem-fin-dim}.  One can easily check that $\Ext
^{i}_{\cat{O}}(\cat{O}_{U},B)=H^{i}(U;B|_{U})$; this is essentially the
definition of sheaf cohomology.  In particular, $S$ has \fhgd \ if and
only if each $\cat{O}_{U}$ has finite projective dimension.  Also, since
the restriction of a flasque sheaf is still flasque and flasque sheaves
have no cohomology, $\Ext ^{i}(\cat{O}_{U},X_{n})=0$ if $i>0$ and $X$ is
a complex of flasque sheaves.  So Lemma~\ref{lem-fin-dim} applies, and
Theorem~\ref{thm-alt-model} gives us the desired model structure.

The characterization of fibrations in Proposition~\ref{prop-fibs}
translates into surjections of \emph{presheaves} with dimensionwise
flasque kernel.  However, sheaf surjections with flasque kernel are also
presheaf surjections, so we get the claimed characterization of
fibrations.
\end{proof}

The author knows of two cases when ringed spaces are guaranteed to have
\fhgd .

\begin{proposition}\label{prop-dim}
Suppose $(S,\cat{O})$ is a ringed space.  
\begin{enumerate}
\item If $S$ is a finite-dimensional noetherian space, then
$(S,\cat{O})$ has \fhgd .
\item If $S$ is a finite-dimensional locally compact topological
manifold that is countable at infinity, in particular if $S$ is a
finite-dimensional compact manifold, then $(S,\cat{O})$ has \fhgd .
\end{enumerate}
\end{proposition}

\begin{proof}
Part~1 follows from the vanishing
theorem~\cite[Theorem~III.2.7]{hartshorne} of Groth\-en\-dieck, since an
open subspace of a finite-dimensional noetherian space is still a
finite-dimensional noetherian space.  Part~2 is an immediate consequence
of~\cite[Proposition~3.2.2]{kashiwara-schapira}.
\end{proof}

We now discuss the cofibrations in the flat model structure.  Recall
that the category $\cat{O}\Mod $ is a closed symmetric monoidal
category.  The monoidal structure is given by the tensor product
$X\otimes _{\cat{O}}Y$, which we will always denote by $X\otimes Y$.
This is defined by forming the obvious presheaf tensor product, and
sheafifying.  On each stalk, the tensor product is the ordinary tensor
product of modules.  In particular, a sheaf $F$ is flat if and only if
each stalk $F_{x}$ is flat as a $\cat{O}_{x}$-module; hence, the sheaves
$\cat{O}_{U}$ are flat.  The closed structure is given by the sheaf Hom;
$\Hom (X,Y)(U)=\cat{O}|_{U}\Mod (X|_{U},Y|_{U})$.  These structures
extend to complexes in the usual way, making $\Ch (\cat{O}\Mod )$ into a
closed symmetric monoidal category.  This works for any symmetric
monoidal additive category, as described
in~\cite[Section~9.2]{hovey-axiomatic}

\begin{definition}\label{defn-DG-flat}
Suppose $\cat{A}$ is a symmetric monoidal abelian category.  Define a
complex $F\in \Ch (\cat{A})$ to be \emph{DG-flat} if each $F_{n}$
is flat, and, for any acyclic complex $K$, the complex $A\otimes K$ is
also acyclic.  
\end{definition}

\begin{proposition}\label{prop-flat-cofibrant}
Suppose $(S,\cat{O})$ is a ringed space with \fhgd .  Then any
cofibration in the flat model structure is a degreewise split
monomorphism on each stalk, with DG-flat cokernel.  
\end{proposition}

We do not know if the converse to this proposition holds, nor even
whether every DG-flat complex is cofibrant.  

\begin{proof}
The maps of $I$ are all degreewise split monomorphisms on each stalk.
Every cofibration is a retract of a transfinite composition of pushouts
of maps of $I$, by the small object
argument~\cite[Theorem~2.1.14]{hovey-model}.  Since retracts,
transfinite compositions, and pushouts all commute with the operation of
taking stalks and preserve split monomorphisms, every cofibration will
be a degreewise split monomorphism on each stalk.  The cokernel of a
cofibration will of course be cofibrant, so to complete the proof it
suffices to show that every cofibrant object is DG-flat.  

Every cofibrant object $A$ is a retract of the colimit of a transfinite
sequence $X_{\alpha }$, where each map $X_{\alpha
}\xrightarrow{}X_{\alpha +1}$ is a pushout of a map of $I$ and
$X_{0}=0$.  Since colimits commute with tensor products and homology, it
suffices to show that, if $X_{\alpha }$ is DG-flat, so is $X_{\alpha
+1}$.  On each stalk, the maps of $I$ are degreewise split monomorphisms
with degreewise flat cokernel, so the same will be true of $X_{\alpha
}\xrightarrow{}X_{\alpha +1}$.  Thus, if $X_{\alpha }$ is a complex of
flat sheaves, so is $X_{\alpha +1}$.  

Now suppose $K$ is an acyclic complex and $f$ is a map of $I$.  Again,
since the maps of $I$ are degreewise split monomorphisms on each stalk,
the map $f\otimes K$ will still be injective.  Thus the map $X_{\alpha
}\otimes K\xrightarrow{}X_{\alpha +1}\otimes K$ will be injective.  It
therefore suffices to show that $f\otimes K$ is a quasi-isomorphism, by
Corollary~\ref{cor-proper}.  In case $f$ is of the form
$D^{n}\cat{O}_{V}\xrightarrow{}D^{n}\cat{O}_{U}$, both the domain and
codomain of $f$ are contractible.  The same will be true of $f\otimes
K$, so $f\otimes K$ will be a quasi-isomorphism.  In case $f$ is of the
form $S^{n-1}\cat{O}_{U}\xrightarrow{}D^{n}\cat{O}_{U}$, the codomain of
$f\otimes K$ is contractible, so it suffices to show that
$S^{n-1}\cat{O}_{U}\otimes K$ is acyclic.  But, since $\cat{O}_{U}$ is
flat, we have $H_{m}(S^{n-1}\cat{O}_{U}\otimes K)= H_{m-n+1}(K)\otimes
\cat{O}_{U}$, so we are done.
\end{proof}

In particular, it follows that cofibrations are pure monomorphisms, in
the sense that, if $f$ is a cofibration and $K$ is an arbitrary complex,
then $f\otimes K$ is still a monomorphism.  

We now show that the flat model structure is compatible with the tensor
product on $\Ch (\cat{O}\Mod )$.  To do this, we need to recall the
definition of this compatibility.  If $f\mathcolon A\xrightarrow{}B$ and
$g\mathcolon C\xrightarrow{}D$ are maps in a cocomplete closed symmetric
monoidal category, we denote the induced map
\[
(A\otimes D) \amalg_{A\otimes C} (B\otimes C) \xrightarrow{} B\otimes D
\]
by $f\boxprod g$.  In case $\cat{C}$ is also a model category, we say
that $\cat{C}$ is a \emph{symmetric monoidal model category} if,
whenever $f$ and $g$ are cofibrations, so is $f\boxprod g$, and
furthermore, if one of $f$ or $g$ is a trivial cofibration, so is
$f\boxprod g$.  This is the condition needed to ensure that the homotopy
category $\ho \cat{C}$ is again closed symmetric monoidal, as explained
in~\cite[Chapter 4]{hovey-model}.  (Actually one also needs a condition
on the unit, but this condition is unnecessary when the unit is
cofibrant, as it is in the flat model structure).

\begin{theorem}\label{thm-boxprod}
Suppose $(S,\cat{O})$ is a ringed space with \fhgd, and $f$ and $g$ are
maps in $\Ch (\cat{O}\Mod )$. 
\begin{enumerate}
\item [(a)] If $f$ is a flat cofibration and $g$ is a monomorphism, then
$f\boxprod g$ is a monomorphism.  
\item [(b)] If $f$ and $g$ are flat cofibrations, then so is $f\boxprod
g$. 
\item [(c)] If $f$ is a flat cofibration and $g$ is an injective
quasi-isomorphism, then $f\boxprod g$ is a quasi-isomorphism.  
\item [(d)] If $f$ is a flat trivial cofibration and $g$ is a
monomorphism, then $f\boxprod g$ is a quasi-isomorphism.  
\end{enumerate}
\end{theorem}

\begin{proof}
As explained in~\cite[Chapter~4]{hovey-model}, since monomorphisms and
injective quasi-isomorphisms are closed under retracts, transfinite
compositions, and pushouts, it suffices to check this
theorem when the flat cofibration is in fact a map of $I$, and the flat
trivial cofibration is a map of $J$.  We begin with parts~(a) and~(c).
Suppose that $g\mathcolon A\xrightarrow{}B$ is a monomorphism, and
suppose $f$ is the map $D^{n}\cat{O}_{V}\xrightarrow{}D^{n}\cat{O}_{U}$.
Let $P$ denote the domain of $f\boxprod g$, and suppose $x\in S$.  If
$x\in V$, then the stalk of $P_{m}$ at $x$ is $(B_{m-n}\oplus
B_{m-n+1})_{x}$; if $x\in U\setminus V$, then the stalk of $P$ at $x$ is
$(A_{m-n}\oplus A_{m-n+1})_{x}$; and if $x$ is not in $U$, then the
stalk of $P$ at $x$ is $0$.  The stalk of the codomain of $f\boxprod g$
at $x$ is $(B_{m-n}\oplus B_{m-n+1})_{x}$ if $x$ is in $U$, and $0$
otherwise, and the map $f\boxprod g$ does the obvious thing on the
stalks.  Hence $f\boxprod g$ is a monomorphism.  Furthermore, the domain
and codomain of $f$ are contractible, so the same will be true for
$f\boxprod g$.  Thus $f\boxprod g$ will be a quasi-isomorphism,
completing the proof of part~(c).  

To complete the proof of part~(a), we must show that $f\boxprod g$ is a
monomorphism, where now $f$ is the map
$S^{n-1}\cat{O}_{U}\xrightarrow{}D^{n}\cat{O}_{U}$.  In this case, the
stalk of the domain $P$ of $f\boxprod g$ at a point $x$ is $0$ if $x\not
\in U$, and otherwise is $(A_{m-n}\oplus B_{m-n+1})_{x}$.  The stalk of
the codomain of $f\boxprod g$ at $x$ is $0$ if $x\not \in U$, and
otherwise is $(B_{m-n}\oplus B_{m-n+1})_{x}$.  The map $f\boxprod g$
does the obvious thing, and so is a monomorphism.  

For part~(c), we can assume $f$ is the map
$S^{n-1}\cat{O}_{U}\xrightarrow{}D^{n}\cat{O}_{U}$.  Then the codomain
of $f\boxprod g$ is contractible, so it suffices to show that the domain
$P$ of $f\boxprod g$ has no homology.  Since $g\mathcolon
A\xrightarrow{}B$ is an injective quasi-isomorphism, and $\cat{O}_{U}$
is flat, $g\otimes S^{n-1}\cat{O}_{U}$ is also an injective
quasi-isomorphism.  Hence its pushout $A\otimes
D^{n}\cat{O}_{U}\xrightarrow{}P$ is also an injective
quasi-isomorphism.  Since $D^{n}\cat{O}_{U}$ is contractible, it follows
that $P$ has no homology.  

Finally, for part~(d), we can assume that both $f$ and $g$ are maps of
$I$.  To calculate $f\boxprod g$ in this case, use the easily checked
(on stalks) fact that $\cat{O}_{U}\otimes \cat{O}_{V}\cong
\cat{O}_{U\cap V}$.  It follows that 
\[
S^{m}\cat{O}_{U}\otimes S^{n}\cat{O}_{V}\cong S^{m+n}\cat{O}_{U\cap V}
\]
and 
\[
S^{m}\cat{O}_{U}\otimes D^{n}\cat{O}_{V}\cong D^{m+n}\cat{O}_{U\cap V}
\]
and that $D^{m}\cat{O}_{U}\otimes D^{n}\cat{O}_{V}$ is an
amalgamation of $D^{m+n-1}\cat{O}_{U\cap V}$ and $D^{m+n}\cat{O}_{U\cap
V}$.  With these identities in hand, the proof is a 
calculation we leave to the reader.  
\end{proof}

\begin{corollary}\label{cor-flat-monoidal}
Suppose $(S,\cat{O})$ is a ringed space with \fhgd .  Then the flat
model structure makes $\Ch (\cat{O}\Mod )$ into a symmetric monoi\-dal
model category.  Furthermore, if $A$ is cofibrant, then the functor
$A\otimes -$ preserves quasi-isomorphisms.  Therefore, to calculate the
derived tensor product up to isomorphism, it suffices to replace
\emph{one} of the factors by a cofibrant complex quasi-isomorphic to it.
\end{corollary}

\begin{proof}
The fact that $\Ch (\cat{O}\Mod )$ is a symmetric monoidal model
category is immediate from Theorem~\ref{thm-boxprod}.  Suppose $A$ is
cofibrant.  Then $A\otimes -$ preserves trivial cofibrations in any
symmetric monoidal model category.  So, in order to see that $A\otimes
-$ preserves quasi-isomorphisms, it suffices to show that, if $p$ is a
trivial fibration, then $A\otimes p$ is a quasi-isomorphism.  Let $K$
denote the kernel of $p$, so that $K$ is an acyclic complex.  Since
cofibrant objects are degreewise flat, $A\otimes p$ is still surjective,
and its kernel is $A\otimes K$.  Since cofibrant objects are DG-flat,
$A\otimes K$ is still acyclic, so the long exact sequence completes the
proof.  In general, the total derived functor of the tensor product is
defined by $X\otimes ^{L}Y = QX \otimes QY$, where $QX$ (resp. $QY$) is
a functorial cofibrant replacement for $X$ (resp. $Y$).  But, since the
map $QX \otimes QY \xrightarrow{}QX\otimes Y$ is a quasi-isomorphism,
$X\otimes ^{L}Y$ is isomorphic in the derived category to $QX\otimes
Y$.  
\end{proof}

Note that Theorem~\ref{thm-boxprod} actually says not only that the flat
model structure is symmetric monoidal, but also that the injective model
structure is a module over the flat model structure, in the sense
of~\cite[Chapter~4]{hovey-model}.  

We can also use Theorem~\ref{thm-boxprod} to conclude that the derived
category of $\cat{O}$-modules is almost a unital algebraic stable homotopy
category~\cite{hovey-axiomatic}.  

\begin{corollary}\label{cor-derived}
Suppose $(S,\cat{O})$ is a ringed space such that $S$ is
a finite-dimensional noetherian space.  Then the derived
category of $\cat{O}$-modules is a symmetric monoidal triangulated
category and $\{\cat{O}_{U} \}$ is a set of small weak generators.  
\end{corollary}

\begin{proof}
It is well-known that the derived category of any abelian category is
triangulated, but this also follows, in a stronger sense of the word
triangulated, from the results of~\cite[Chapter~7]{hovey-model}.  We
have already seen that the flat model caetgory is a symmetric monoidal
model category, so the derived category is also closed symmetric
monoidal in a way that is compatible with the triangulation
(see~\cite[Chapter~6]{hovey-model}, with one technical point dealt with
by~\cite[Corollary~5.6.10]{hovey-model}).  Since the flat model
structure is cofibrantly generated, the cofibers of the generating
cofibrations form a set of weak
generators~\cite[Section~7.3]{hovey-model}.  In our case, these are the
objects $S^{n}\cat{O}_{U}$ (the cofibers of the maps of $J$ are trivial
in the derived category).  Because $S$ is noetherian, the presheaf
colimit of a direct system of sheaves coincides with the sheaf
colimit~\cite[Exercise~II.1.11]{hartshorne}.  It follows from this that
$\Ch (\cat{O}\Mod )(S^{n}\cat{O}_{U},-)$ commutes with direct colimits.
The results of~\cite[Section~7.4]{hovey-model} then show that
$S^{n}\cat{O}_{U}$ is small (in the triangulated sense) in the derived
category.
\end{proof}

The derived category of $\cat{O}$-modules is known to be a symmetric
monoidal triangulated category even without the \fhgd \
assumption~\cite[Section~2.5]{lipman}.  This might indicate that there
is some replacement for the flat model structure that works more
generally, or it might indicate that model categories are simply not
adequate to cope with the general case.

To show that the derived category is in fact a unital algebraic stable
homotopy category, we would need to know that the generators
$\cat{O}_{U}$ are strongly dualizable.  This would mean we would need to
show that the natural map 
\[
R\Hom (\cat{O}_{U},\cat{O})\otimes \cat{O}_{V}\xrightarrow{}R\Hom
(\cat{O}_{U},\cat{O}_{V})
\]
is a quasi-isomorphism.  (We don't need the derived tensor product since
$\cat{O}_{V}$ is cofibrant, using Corollary~\ref{cor-flat-monoidal}).
Unfortunately, this is false in the simplest non-trivial example.
Indeed, consider sheaves of abelian groups on the Sierpinski space $S$.
Recall $S$ has only two points, exactly one of which is open.  Sheaves
on $S$ coincide with presheaves, which in turn coincide with maps
$A\xrightarrow{}B$ of abelian groups.  In this case, $\cat{O}$ is the
identity map $\Z \xrightarrow{}\Z $, and, taking $U$ to be the open
point, $\cat{O}_{U}$ is the map $0\xrightarrow{}\Z $.  One can then
calculate to find $R\Hom (\cat{O}_{U},\cat{O})=R\Hom
(\cat{O}_{U},\cat{O}_{U})$, but this equality is destroyed on tensoring
the left hand side with $\cat{O}_{U}$.   A similar counterexample works
if we think of $S$ as the underlying space of $\Spec \Z _{(p)}$. 

There is an additional condition that a symmetric monoidal model
category might satisfy, called the \emph{monoid
axiom}~\cite{schwede-shipley-monoids}.  This axiom guarantees that the
monoids in a symmetric monoidal model category, and the modules over a
given monoid, themselves form model categories.  The monoid axiom
asserts that every map in $K\cof $ is a weak equivalence, where $K$ is
the class consisting of all maps $f\otimes X$, where $f$ is a trivial
cofibration and $X$ is an arbitrary object.  

\begin{theorem}\label{thm-monoid}
Suppose $(S,\cat{O})$ is a ringed space with \fhgd .  The the flat model
structure on $\Ch (\cat{O}\Mod )$ satisfies the monoid axiom.
\end{theorem}

\begin{proof}
Suppose $f$ is a flat trivial cofibration, and $X$ is an arbitrary
object.  Then $0\xrightarrow{}X$ is a monomorphism, so applying
Theorem~\ref{thm-boxprod} shows that $f\otimes X$ is an injective
quasi-isomorphism.  Corollary~\ref{cor-homology-J} completes the proof.
\end{proof}

The following corollary follows immediately from
Theorem~\ref{thm-monoid} and~\cite{schwede-shipley-monoids}.  

\begin{corollary}\label{cor-monoids}
Suppose $(S,\cat{O})$ is a ringed space with \fhgd .  Then:
\begin{enumerate}
\item [(a)] The category of monoids in $\Ch (\cat{O}\Mod )$ is a
cofibrantly generated model category, where a map of monoids is a weak
equivalence or a fibration if and only if it is so in the flat model
structure on $\Ch (\cat{O}\Mod )$. 
\item [(b)] Given a monoid $R$ in $\Ch (\cat{O}\Mod )$, the category of
$R$-modules, $R\Mod $, is a cofibrantly generated proper model category,
where a map of modules is a weak equivalence or a fibration if and only
if it is so in the flat model structure on $\Ch (\cat{O}\Mod )$.  
\item [(c)] If $R$ is a commutative monoid, then $R\Mod $ is a symmetric
monoidal model category satisfying the monoid axiom.  Furthermore, the
category of algebras over $R$ is a cofibrantly generated model category,
where a map of algebras is a weak equivalence or fibration if and only
if it is so in the flat model structure on $\Ch (\cat{O}\Mod )$.  
\end{enumerate}
\end{corollary}

The category of monoids will certainly be right proper, but we do not
know if the category of monoids is left proper.  

A map of monoids $R\xrightarrow{}R'$ will induce the usual induction and
restriction adjunction $R\Mod \xrightarrow{}R'\Mod $.  This adjunction
will be a Quillen adjunction, but we would like something more. 

\begin{proposition}\label{prop-weak-equiv-of-monoids}
Suppose $(S,\cat{O})$ is a ringed space with \fhgd , and
$R\xrightarrow{}R'$ is a weak equivalence of monoids in $\Ch
(\cat{O}\Mod )$.  Then the induction and restriction adjunction $R\Mod
\xrightarrow{}R'\Mod $ is a Quillen equivalence.
\end{proposition}

\begin{proof}
It suffices to show that, if $N$ is a cofibrant $R$-module, then
$-\otimes _{R}N$ preserves weak equivalences,
by~\cite{schwede-shipley-monoids}.  The proof of this is very similar to
the proof of the corresponding fact in
Corollary~\ref{cor-flat-monoidal}, so we leave it to the reader.
\end{proof}

We now investigate the functoriality of the flat model structure.
Suppose we have a map of ringed spaces $f\mathcolon
(S,\cat{O}_{S})\xrightarrow{}(T,\cat{O}_{T})$.  Recall that this is a
continuous map $f\mathcolon S\xrightarrow{}T$ together with a map of
sheaves of rings $\cat{O}_{T}\xrightarrow{}f_{*}\cat{O}_{S}$.  Here, for
any sheaf $X$ on $S$, $f_{*}(X)$ is the sheaf on $T$ defined by
$f_{*}(X)(U)=X(f^{-1}(U))$.  If $X$ is an $\cat{O}_{S}$-module, then
$f_{*}X$ is an $f_{*}\cat{O}_{S}$-module, and so an $\cat{O}_{T}$-module
by restriction.  The functor $f_{*}\mathcolon \cat{O}_{S}\Mod
\xrightarrow{}\cat{O}_{T}\Mod $ has a left adjoint $f^{*}$.  To define
this, recall that if $Y$ is a sheaf on $T$, $f^{-1}Y$ is the sheaf on
$S$ associated to the presheaf that takes $U$ to $\colim _{V\supseteq
f(U)}Y(V)$.  The functor $f^{-1}$ is left adjoint to $f_{*}$ on the
category of sheaves of abelian groups, so in particular we have a map
$f^{-1}\cat{O}_{T}\xrightarrow{}\cat{O}_{S}$.  Given an
$\cat{O}_{T}$-module $Y$, we define $f^{*}Y = \cat{O}_{S}\otimes
_{f^{-1}\cat{O}_{T}} f^{-1}Y$.  It is well-known that $f^{*}$ is left
adjoint to $f_{*}$ and is symmetric
monoidal~\cite[Section~0.4.3]{grothendieck-ega1}.

One can verify using adjointness that, if $U$ is an open subset of $T$, then
$f^{*}\cat{O}_{U}=\cat{O}_{f^{-1}U}$.  Hence we have the following
proposition.  

\begin{proposition}\label{prop-flat-functor}
Suppose $f\mathcolon (S,\cat{O}_{S})\xrightarrow{}(T,\cat{O}_{T})$ is a
map of ringed spaces with \fhgd.  Then $f^{*}$ is a left Quillen functor
with respect to the flat model structures.
\end{proposition}

In particular, this shows that the total left derived functor of $f^{*}$
exists and is left adjoint to the total right derived functor of
$f_{*}$.  It is proved in~\cite[Section~2.7]{lipman} that the total left
derived functor of $f^{*}$ exists without the \fhgd \  hypotheses.  It is
disconcerting that we are unable to reproduce this result using model
categories.  


\providecommand{\bysame}{\leavevmode\hbox to3em{\hrulefill}\thinspace}

\end{document}